# GLOBAL STABILITY RESULTS FOR TRAFFIC NETWORKS

**Iasson Karafyllis[*] and Markos Papageorgiou[**]**

**[*]Dept. of Mathematics, National Technical University of Athens, Zografou Campus, 15780, Athens, Greece (email: iasonkar@central.ntua.gr )**

**[**]School of Production Engineering and Management, Technical University of Crete, Chania, 73100, Greece (email: markos@dssl.tuc.gr )**

**Abstract**

This paper provides sufficient conditions for global asymptotic stability and global exponential stability, which can be applied to nonlinear, large-scale, uncertain discrete-time systems. The conditions are derived by means of vector Lyapunov functions. The obtained results are applied to traffic networks for the derivation of sufficient conditions of global exponential stability of the uncongested equilibrium point of the network. Specific results and algorithms are provided for freeway models. Various examples illustrate the applicability of the obtained results.

**Keywords:** nonlinear systems, discrete-time systems, traffic networks.

## 1. Introduction

The purpose of the present paper is three-fold:

- to provide sufficient conditions for Global Asymptotic Stability (GAS) and Global Exponential Stability (GES), which can be easily applied to nonlinear, large-scale, uncertain discrete-time systems;
- to apply the aforementioned sufficient conditions to traffic networks and obtain conditions, which guarantee the GES of the uncongested equilibrium point;
- to study the stability properties of freeway models and obtain easily checkable conditions which guarantee the GES of the uncongested equilibrium point.

Vector Lyapunov functions are particularly useful to large-scale discrete-time systems. Sufficient stability conditions by means of vector Lyapunov functions have been proposed in [9]. More recently, small-gain conditions have been proposed in [17], which can be expressed by means of a vector Lyapunov function formulation (as shown in [11]). In this work, we propose a set of conditions expressed by means of vector Lyapunov functions, which guarantee GAS and GES (Theorem 2.3) and can be applied easily to nonlinear, large-scale, uncertain discrete-time systems. The basis for the applicability is the expression of the stability condition by means of a condition on the spectral radius of a nonnegative matrix. Therefore, we can apply recent results on nonnegative matrices that provide upper bounds for the spectral radius (see [3]). The stability notions used in this work are the standard stability notions for discrete-time systems used in [9,10,11,13], but we also allow the discrete-time, uncertain system to be defined on a subset of a finite-dimensional space. Discrete-time systems defined on a subset of a finite-dimensional space were studied in [23].

The conservatism of the obtained stability conditions can be reduced significantly if we have an accurate description of a trapping region of the system: this feature is exploited throughout the present work. A nonlinear system with a trapping region is a system for which all solutions enter a specific set after an initial transient period (for continuous-time systems without inputs the name “global uniform ultimate boundedness” is used in [12] when the corresponding set is compact; the term “dissipative system” is used in the literature of continuous-time dynamical systems with compact corresponding sets; see [23] and the discussion on page 22 of the book [24]).

The obtained stability results are applied to traffic networks (Section 3). More specifically, we develop a general model for traffic networks, which consists of an arbitrary number of elementary components. The components can be interconnected to form any two-dimensional structure for the overall traffic network. This general formulation allows for a plethora of diverse traffic network infrastructures to be addressed on the basis of a unifying modeling approach.

In particular, the traffic network structures and problems that can be considered as special cases of the proposed network model include: urban road networks consisting of interconnected links which are modelled as store-and-forward components [1] or cell-transmission links [4]; large urban networks consisting of smaller homogeneous sub-networks which are modelled by use of individual NFD (network fundamental diagrams) [2]; freeway stretches or networks consisting of series of links which are modelled via the discretized LWR (Lighthill-Whitham-Richards) model [14] or its simplified CTM (Cell Transmission Model) version [6]; freeway stretches or networks consisting of appropriately interconnected lane-segments which provide a basis for modelling (partially) automated highway traffic [21]; large mixed (corridor) networks consisting of urban and freeway links [19]. As a matter of fact, the same generic approach may also be used for modelling water networks consisting of interconnected links which are modelled by discretized versions of the Lighthill-Whitham model [16], see [5,18]. Our main related result (Theorem 3.1) provides explicit formulas for the elements of a specific nonnegative matrix whose spectral radius is critical for the GES of the uncongested equilibrium point of the traffic network. Therefore, our results can be used in a straightforward way for the determination of the stability properties in a given traffic network within this general framework (see Example 3.3).

The obtained results are next specialized to the case of a freeway stretch (Section 4). As mentioned, the overall model in this constellation consists of a series of subsequent segments (cells) and is very similar to (but more general than) the known first-order discrete Godunov approximations (see [7]) to the kinematic-wave partial differential equation of the LWR-model (see [16,20]) with nonlinear ([14]) or piecewise linear (CTM, [6]) outflow functions. However, the presented framework can also accommodate recent modifications of the LWR-model as in [15] to reflect the so-called capacity drop phenomenon. Our main related result (Corollary 4.4) provides an easily implementable algorithm for the determination of the stability properties of the uncongested equilibrium point of the freeway stretch. The results are different from other results in the literature on the Cell Transmission Model (see [8]), since our model is more general than the Cell Transmission Model and our methodology is completely different from the methodology used in [8].

**Notation.** Throughout this paper, we adopt the following notation:

* $\Re_+ := [0,+\infty)$. For every set $S$, $S^n = \underbrace{S \times \ldots \times S}_{n\ times}$ for every positive integer $n$. $\Re_+^n := (\Re_+)^n$.

* By $C^0(A\,;\Omega)$, we denote the class of continuous functions on $A \subseteq \Re^n$, which take values in $\Omega \subseteq \Re^m$. By $C^k(A\,;\Omega)$, where $k \ge 1$ is an integer, we denote the class of functions on $A \subseteq \Re^n$ with continuous derivatives of order $k$, which take values in $\Omega \subseteq \Re^m$.

* We say that a function $\rho : \Re_+ \to \Re_+$ is positive definite if $\rho(0) = 0$ and $\rho(s) > 0$ for all $s > 0$. By $K$ we denote the set of positive definite, increasing and continuous functions. We say that a positive definite, increasing and continuous function $\rho : \Re_+ \to \Re_+$ is of class $K_\infty$ if $\lim_{s \to +\infty} \rho(s) = +\infty$. By $KL$ we denote the set of all continuous functions $\sigma = \sigma(s,t) : \Re_+ \times \Re_+ \to \Re_+$ with the properties: (i) for each $t \ge 0$ the mapping $\sigma(\cdot\,,t)$ is of class $K$ ; (ii) for each $s \ge 0$, the mapping $\sigma(s,\cdot)$ is non-increasing with $\lim_{t \to +\infty} \sigma(s,t) = 0$.

* Let $x, y \in \Re^n$. We say that $x \le y$ if and only if $(y - x) \in \Re_+^n$.

* Let $x \in \Re^n$. By $|x|$ we denote the Euclidean norm of $x \in \Re^n$. Let $A \in \Re^{n \times n}$ be a real matrix. By $|A|$ we denote the induced matrix norm, i.e., $|A| = \sup\left\{ |Ax| : x \in \Re^n, |x| = 1 \right\}$. The spectral radius $\rho(A)$ of the square matrix $A \in \Re^{n \times n}$ is equal to $\max_{i=1,\ldots,n} |\lambda_i|$, where $\lambda_1,\ldots,\lambda_n$ are the eigenvalues of $A \in \Re^{n \times n}$. When all elements of the matrix $A \in \Re^{n \times n}$ are non-negative, then we say that the matrix $A$ is non-negative and we write $A \in \Re_+^{n \times n}$.

* For every $x \in \Re$, $[x]$ denotes the integer part of $x \in \Re$, i.e., the largest integer $k$ that satisfies $k \le x$.

## 2. Vector Lyapunov Stability Criteria for Discrete-Time Networks

Consider the discrete-time system:

$$\begin{aligned} & x^+ = F(d,x) \\ & x \in S \subseteq \Re^n , d \in D \end{aligned} \tag{2.1}$$

where $S \subseteq \Re^n$ is a non-empty closed set with $x^* \in S$, $D \subseteq \Re^l$ is a non-empty, compact set, $F : D \times S \to S$ is a locally bounded mapping, being continuous on the set $D \times \{x^*\}$ with $F(d,x^*) = x^*$ for all $d \in D$. We suppose that for every $\delta > 0$ it holds that $\left\{ x \in S : 0 < \left| x - x^* \right| \leq \delta \right\} \neq \emptyset$.

In order to develop the Vector Lyapunov Stability criteria we need the notion of a trapping region. A nonlinear system with a trapping region is a system for which all solutions enter a specific set after an initial transient period (for continuous-time systems without inputs the name "global uniform ultimate boundedness" is used in [12] when the corresponding set is compact; the term "dissipative system" is used in the literature of continuous-time dynamical systems with compact corresponding sets; see [23] and the discussion on page 22 of the book [24]).

**Definition 2.1:** *A trapping region for system (2.1) is a set $A \subseteq S$ for which there exists a non-negative integer $m \geq 0$ such that for every $x_0 \in S$, $\{d_i \in D\}_{i=0}^{\infty}$, the solution $x(t)$ of (2.1) with initial condition $x(0) = x_0$ corresponding to input $\{d_i \in D\}_{i=0}^{\infty}$ satisfies $x(t) \in A$ for all $t \geq m$.*

A direct consequence of Definition 2.1 is that every trapping region for (2.1) must contain the equilibrium point $x^* \in S$. We next define rigorously the robust stability notions used for system (2.1).

**Definition 2.2:** *We say that $x^* \in S$ is Robustly Globally Asymptotically Stable (RGAS) for system (2.1), if there exists a function $\sigma \in KL$ such that for every $x_0 \in S$, $\{d_i \in D\}_{i=0}^{\infty}$, the solution $x(t)$ of (2.1) with initial condition $x(0) = x_0$ corresponding to input $\{d_i \in D\}_{i=0}^{\infty}$ satisfies $\left| x(t) - x^* \right| \leq \sigma\left( \left| x_0 - x^* \right|, t \right)$ for all integers $t \geq 0$. We say that $x^* \in S$ is Robustly Globally Exponentially Stable (RGES) for system (2.1) if there exist constants $M, \sigma > 0$ such that for every $x_0 \in S$, $\{d_i \in D\}_{i=0}^{\infty}$, the solution $x(t)$ of (2.1) with initial condition $x(0) = x_0$ corresponding to input $\{d_i \in D\}_{i=0}^{\infty}$ satisfies $\left| x(t) - x^* \right| \leq M \exp(-\sigma t) \left| x_0 - x^* \right|$ for all integers $t \geq 0$.*

When the system is free of disturbances then we omit the adjective "Robust" and we say that an equilibrium point is Globally Asymptotically Stable (GAS) or Globally Exponentially Stable (GES).

We are now ready to state the main result of the section.

**Theorem 2.3:** *Consider system (2.1) and suppose that $A \subseteq S$ is a trapping region for system (2.1). Moreover, suppose that there exists a family of functions $V_i : S \to \Re_+$ $(i = 1,...,l)$, a matrix $\Gamma = \left\{ \gamma_{i,j} \geq 0, i,j = 1,...,l \right\} \in \Re_+^{l \times l}$ and functions $a_1, a_2 \in K_\infty$ with $a_1(s) \leq a_2(s)$ for all $s \geq 0$ such that the following inequalities hold for all $x \in A, d \in D$ and $i = 1,...,l$:*

$$a_1\left( \left| x - x^* \right| \right) \leq \max_{i=1,...,l} \left( V_i(x) \right) \leq a_2\left( \left| x - x^* \right| \right) \tag{2.2}$$

$$V_i(F(d,x)) \leq \sum_{j=1}^{l} \gamma_{i,j} V_j(x) \tag{2.3}$$

*Moreover, suppose that the spectral radius $\rho(\Gamma)$ of the matrix $\Gamma = \left\{ \gamma_{i,j} \geq 0, i,j = 1,...,l \right\} \in \Re_+^{l \times l}$ is less than 1. Then $x^* \in S$ is RGAS for system (2.1). Moreover, if there exist constants $L \geq 0$, $0 < K_1 \leq K_2$, $p > 0$ such that $\max\left\{ \left| F(d,x) - x^* \right| : d \in D \right\} \leq L \left| x - x^* \right|$ for all $x \in S \setminus A$ and if $a_i(s) = K_i s^p$ $(i = 1,2)$ for all $s \geq 0$ then $x^* \in S$ is RGES for system (2.1).*

Since the matrix $\Gamma = \{\gamma_{i,j} \ge 0, i,j = 1,...,l\} \in \Re_+^{l \times l}$ involved in the statement of the Theorem 2.3 is non-negative, there are effective computational tools, which can guarantee that its spectral radius is less than 1 (see [3]). For example, if

$$\max_{i=1,...,n} \left( \sum_{j=1}^{n} \gamma_{i,j} \right) < 1$$

or if there exists $\varepsilon > 0$ such that

$$\max_{i=1,...,n} \left( \frac{\sum_{j=1}^{n} (\varepsilon + \gamma_{i,j}) \sum_{k=1}^{n} (\varepsilon + \gamma_{j,k})}{n\varepsilon + \sum_{j=1}^{n} \gamma_{i,j}} \right) < 1$$

then the spectral radius of $\Gamma = \{\gamma_{i,j} \ge 0, i,j = 1,...,l\} \in \Re_+^{l \times l}$ is less than 1. The above conditions can be used for large-scale systems easily.

**Proof:** Let $x_0 \in S$, $\{d_i \in D\}_{i=0}^{\infty}$ (arbitrary) be given and consider the solution $x(t)$ of (2.1) with initial condition $x(0) = x_0$ corresponding to input $\{d_i \in D\}_{i=0}^{\infty}$. Let $m \ge 0$ be the non-negative integer involved in Definition 2.1. Let $j \in \{0,...,m\}$ be the smallest non-negative integer for which it holds that $x(t) \in A$ for all $t \ge m$ (the fact that $j$ exists and satisfies $j \in \{0,...,m\}$ is a direct consequence of the fact that $A \subseteq S$ is a trapping region for system (2.1).

We next show that there exists a function $b \in K_\infty$ with $b(s) \ge 0$ for all $s \ge 0$, such that

$$\max_{k=0,...,j} \left| x(k) - x^* \right| \le b\left( \left| x_0 - x^* \right| \right) \tag{2.4}$$

Indeed, if there exists a constant $L \ge 0$ such that $\max\left\{ \left| F(d,x) - x^* \right| : d \in D \right\} \le L \left| x - x^* \right|$ for all $x \in S \setminus A$, then we may define $b(s) := \max(1, L^m)s$ for all $s \ge 0$. The fact that (2.4) holds is a direct consequence of the fact that $j \le m$, equation (2.1) and the resulting inequality $\left| x(t+1) - x^* \right| \le \max(1, L) \left| x(t) - x^* \right|$ which holds for all $t = 0,...,j-1$, for the case that $j \ge 1$. When $j = 0$, then (2.4) holds automatically.

For the general case, we define:

$$a(s) := \sup\left\{ \left| F(d,x) - x^* \right| : (d,x) \in D \times S, \left| x - x^* \right| \le s \right\}. \tag{2.5}$$

Clearly, $a(s)$ is well-defined by (2.5) for all $s \ge 0$, since $F : D \times S \to S$ is a locally bounded mapping. Continuity of $F : D \times S \to S$ on the set $D \times \{x^*\}$ in conjunction with the fact that $D \subseteq \Re^l$ is a compact set with $F(d, x^*) = x^*$ for all $d \in D$ implies that

$$\lim_{s \to 0^+} a(s) = a(0) = 0 . \tag{2.6}$$

The function $a : \Re_+ \to \Re_+$ is non-decreasing. We next define the function $\tilde{a} : \Re_+ \to \Re_+$ :

$$\tilde{a}(s) := \frac{1}{s} \int_{s}^{2s} a(w) dw, \text{ for } s > 0 \text{ and } \tilde{a}(0) := 0 . \tag{2.7}$$

The integral appearing in (2.7) is the Riemann integral of $a$ (well-defined because $a : \Re_+ \to \Re_+$ is non-decreasing). Monotonicity of $a : \Re_+ \to \Re_+$ and definition (2.7) implies that $a(s) \le \tilde{a}(s) \le a(2s)$ for all $s \ge 0$. Consequently, definition (2.7) in conjunction with (2.5), (2.6) and the previous inequality shows that $\tilde{a}:\Re_+ \to \Re_+$ is a non-decreasing, continuous function which satisfies:

$$\max\left\{\left|F(d,x)-x^*\right| : d \in D\right\} \le \tilde{a}\left(\left|x-x^*\right|\right) \text{ for all } x \in S \tag{2.8}$$

Next, we define:

$$b := \underbrace{\bar{a} \circ \ldots \circ \bar{a}}_{m \text{ times}}, \text{ where } \bar{a}(s) := \tilde{a}(s) + s \text{ for all } s \ge 0 \tag{2.9}$$

Definition (2.9) shows that $b \in K_\infty$ with $b(s) \ge 0$ for all $s \ge 0$. Using the fact that $j \le m$, inequality $\left|x(t+1)-x^*\right| \le \tilde{a}\left(\left|x(t)-x^*\right|\right)$ which holds for all $t = 0,...,j-1$ (a direct consequence of (2.8) and (2.1)) for the case that $j \ge 1$ and definition (2.9), we obtain (2.4). When $j = 0$, then (2.4) holds automatically.

Since the spectral radius $\rho(\Gamma)$ of the matrix $\Gamma = \left\{\gamma_{i,j} \ge 0, i,j = 1,...,l\right\} \in \Re_+^{l \times l}$ is less than 1, it follows that there exist constants $M \ge 1$, $\sigma > 0$ such that

$$\left|\Gamma^t\right| \le M \exp(-\sigma t), \text{ for all integers } t \ge 0 \tag{2.10}$$

(see [22]). Next define:

$$\xi(t) = (V_1(x(t+j)),...,V_l(x(t+j)))' \in \Re_+^l \text{ for all } t \ge 0 \tag{2.11}$$

Equation (2.1) in conjunction with inequalities (2.3) imply that the following recursive relation holds for all $t \ge 0$:

$$\xi(t+1) \le \Gamma \xi(t) \tag{2.12}$$

The inequality in (2.12) is the usual partial order of $\Re^l$ (i.e., $x \le y$ for two vectors $x, y \in \Re^l$ iff $(y-x) \in \Re_+^l$). Using the fact that $\Gamma = \left\{\gamma_{i,j} \ge 0, i,j = 1,...,l\right\} \in \Re_+^{l \times l}$ is a non-negative matrix (and consequently a monotone operator which satisfies $\Gamma x \le \Gamma y$ for all vectors $x, y \in \Re^l$ with $x \le y$), we obtain from (2.12):

$$\xi(t) \le \Gamma^t \xi(0), \text{ for all } t \ge 0 \tag{2.13}$$

Using (2.10), (2.13), definition (2.11) and (2.2), we get:

$$a_1\left(\left|x(j+t)-x^*\right|\right) \le M \exp(-\sigma t)\sqrt{l}\, a_2\left(\left|x(j)-x^*\right|\right), \text{ for all } t \ge 0 \tag{2.14}$$

Using (2.4) and (2.14), we get:

$$a_1\left(\left|x(j+t)-x^*\right|\right) \le M \exp(-\sigma t)\sqrt{l}\, a_2\left(b\left(\left|x_0-x^*\right|\right)\right), \text{ for all } t \ge 0 \tag{2.15}$$

Since $a_1(s) \le a_2(s)$ for all $s \ge 0$, and since $M \ge 1$, $j \le m$, it follows from (2.4), (2.15) that the following estimate holds for all $t \ge 0$:

$$a_1\left(\left|x(t)-x^*\right|\right) \le M \exp(-\sigma(t-m))\sqrt{l}\, a_2\left(b\left(\left|x_0-x^*\right|\right)\right) \tag{2.16}$$

Inequality (2.16) shows that the estimate $\left|x(t)-x^*\right| \le \sigma\left(\left|x_0-x^*\right|, t\right)$ holds for all $t \ge 0$ with $\sigma(s,t) := a_1^{-1}\left(M \exp(-\sigma(t-m))\sqrt{l}\, a_2(b(s))\right)$ (notice that $\sigma \in KL$) and consequently $x^* \in S$ is RGAS for system (2.1).

If there exist constants $L \ge 0$, $0 < K_1 \le K_2$, $p > 0$ such that $\max\left\{\left|F(d,x)-x^*\right| : d \in D\right\} \le L\left|x-x^*\right|$ for all $x \in S \setminus A$ and if $a_i(s) = K_i s^p$ ($i = 1,2$) for all $s \ge 0$ then inequality (2.15) implies that

$$\left|x(j+t)-x^*\right| \le \left(\frac{M\sqrt{l}\,K_2}{K_1}\right)^{1/p} \exp(-\sigma t/p)\max\left(1,L^m\right)\left|x_0-x^*\right|, \text{ for all } t \ge 0 \tag{2.17}$$

The reader should notice that we have used the fact that in this case (2.4) holds with $b(s):=\max(1,L^m)s$ for all $s \ge 0$. Again, it follows from (2.17), (2.4) with $b(s):=\max(1,L^m)s$ for all $s \ge 0$ and the facts that $j \le m$, $0 < K_1 \le K_2$ that the following estimate holds for all $t \ge 0$:

$$\left|x(t)-x^*\right| \le \left(\frac{M\sqrt{l}\,K_2}{K_1}\right)^{1/p} \exp(-\sigma(t-m)/p)\max\left(1,L^m\right)\left|x_0-x^*\right|, \text{ for all } t \ge 0 \tag{2.18}$$

which directly implies that $x^* \in S$ is RGES for system (2.1). The proof is complete. $\lhd$

**Example 2.4:** Consider the large-scale discrete-time system

$$x_1^+ = (1-2r)x_1 + dx_2 \tag{2.19}$$

$$x_i^+ = dx_{i-1} + (1-2r)x_i + dx_{i+1}, \ i = 2,...,n-1 \tag{2.20}$$

$$x_n^+ = dx_{n-1} + (1-2r)x_n \tag{2.21}$$

where $d \in D := [-r, r]$ and $r > 0$ is a constant. Systems of the form (2.19), (2.20), (2.21) with $d \equiv r$ arise when a finite difference numerical scheme are applied to the heat equation $\frac{\partial x}{\partial t}(t,z) = D\frac{\partial^2 x}{\partial z^2}(t,z)$ on $t > 0$, $z \in (0,L)$ with Dirichlet boundary conditions $x(t,0) = x(t,L) = 0$, where $D, L > 0$ are constants. In this case, the parameter $r > 0$ is equal to $r = \frac{D\delta t}{(\delta z)^2}$, where $\delta t$ is the discretization time step and $\delta z = \frac{L}{n+1}$ is the discretization space step.

System (2.19), (2.20), (2.21) is a system of the form (2.1) with $S = \Re^n$, $x^* = 0 \in \Re^n$. The stability properties of system (2.19), (2.20), (2.21) can be studied by means of Theorem 2.3, the trapping region $A = S = \Re^n$ and the family of functions $V_i(x) = |x_i|$ for $i = 1,...,n$. The matrix $\Gamma = \{\gamma_{i,j} : i, j = 1,...,n\}$ is the tridiagonal matrix defined by means of the formulae:

$$\gamma_{i,i} := |1-2r|, \text{ for } i = 1,...,n \tag{2.22}$$

$$\gamma_{i,i+1} := r \text{ for } i = 1,...,n-1 \tag{2.23}$$

$$\gamma_{i,i-1} := r, \text{ for } i = 2,...,n. \tag{2.24}$$

When $r \le \frac{1}{2}$ the eigenvalues of $\Gamma \in \Re_+^{n\times n}$ are $\lambda_i = 1-2r+2r\cos\left(\frac{i\pi}{n+1}\right)$, $i = 1,...,n$. Therefore, for $r \le \frac{1}{2}$ the spectral radius $\rho(\Gamma)$ of $\Gamma \in \Re_+^{n\times n}$ is less than 1. Since (2.2) holds with $a_1(s) := \frac{s}{\sqrt{n}}$ and $a_2(s) := s$ for all $s \ge 0$, we conclude from Theorem 2.3 that $x^* = 0 \in \Re^n$ is RGES. The condition $r \le \frac{1}{2}$ is exactly the Courant-Friedrichs-Lewy condition for the stability of the finite difference numerical scheme for the heat equation (see the discussion on pages 206-208 in [13]). $\lhd$

## 3. Global Stability Results for Traffic Networks

This section is devoted to the derivation of sufficient conditions that guarantee RGES for the equilibrium point of a traffic network. We consider a traffic network which consists of $n$ components. The number of vehicles at time $t \geq 0$ in component $i \in \{1,...,n\}$ is denoted by $x_i(t)$. The outflow and the inflow of vehicles of the component $i \in \{1,...,n\}$ at time $t \geq 0$ are denoted by $q_i(t) \geq 0$ and $F_i(t) \geq 0$, respectively. Consequently, the balance of vehicles for each component $i \in \{1,...,n\}$ gives:

$$x_i(t+1) = x_i(t) - q_i(t) + F_i(t)\,,\ i = 1,...,n\,,\ t \geq 0\,. \tag{3.1}$$

Each component of the network has capacity $a_i > 0$ ($i = 1,...,n$). Our first assumption states that the inflow of vehicles at the component $i \in \{1,...,n\}$ at time $t \geq 0$, denoted by $F_i(t) \geq 0$, cannot exceed the number of free positions for vehicles of component $i \in \{1,...,n\}$ at time $t \geq 0$, i.e.,

$$F_i(t) = \min\left(a_i - x_i(t), \widetilde{F}_i(t)\right),\ i = 1,...,n\,,\ t \geq 0 \tag{3.2}$$

where $\widetilde{F}_i(t) \geq 0$ is the attempted inflow of vehicles at the component $i \in \{1,...,n\}$ at time $t \geq 0$.

Our second assumption is dealing with the attempted outflows and inflows. We assume that there exist functions $f_i : D \times [0, a_i] \to \Re_+$ with $f_i(d, x_i) \leq x_i$ for all $(d, x_i) \in D \times [0, a_i]$, where $D \subseteq \Re^l$ is a non-empty, compact set, non-negative constants $p_{i,j} \geq 0$, $i, j = 1,...,n$ with $p_{i,i} = 0$ for $i = 1,...,n$ and non-negative constants $Q_i \geq 0$, $i = 1,...,n$ so that:

$$\begin{pmatrix} \textit{attempted flow of vehicles} \\ \textit{from component i to component j} \end{pmatrix} = p_{i,j} f_i(d, x_i)\,,\ i, j = 1,...,n \tag{3.3}$$

$$\begin{pmatrix} \textit{attempted flow of vehicles} \\ \textit{from component i to regions out of the network} \end{pmatrix} = Q_i f_i(d, x_i)\,,\ i = 1,...,n \tag{3.4}$$

We also assume that:

$$\sum_{j=1}^{n} p_{i,j} + Q_i = 1\,. \tag{3.5}$$

In other words the total attempted outflow from component $i \in \{1,...,n\}$ is equal to $f_i(d, x_i)$. Equations (3.3) (or (3.4)) imply that the attempted flows from component $i \in \{1,...,n\}$ of the network to another component of the network (or to regions out of the network) are constant proportions of the total attempted outflow from component $i \in \{1,...,n\}$. In traffic terminology, $p_{i,j}$ are turning rates at network-internal bifurcations, while $Q_i$ are exit rates.

Let $v_i > 0$ ($i = 1,...,n$) denote the attempted inflow to component $i \in \{1,...,n\}$ from the region out of the network. Our assumptions lead us to the following equations:

$$\widetilde{F}_i(t) = v_i + \sum_{j=1}^{n} p_{j,i} f_j(d(t), x_j(t))\,,\ i = 1,...,n\,,\ t \geq 0\,. \tag{3.6}$$

Equations (3.2) and (3.6) imply that the percentage of the attempted inflow of vehicles at the component $i \in \{1,...,n\}$ at time $t \geq 0$, which becomes actual inflow of vehicles at the component $i \in \{1,...,n\}$ at time $t \geq 0$, denoted by $s_i(t) \in [0,1]$ is given by:

$$s_i(t)=\frac{\min\left(a_i-x_i(t),v_i+\sum_{j=1}^{n}p_{j,i}f_j(d(t),x_j(t))\right)}{v_i+\sum_{j=1}^{n}p_{j,i}f_j(d(t),x_j(t))}\,,\; i=1,...,n\,,\; t\geq 0\,. \tag{3.7}$$

Our final assumption relates the actual inflows with the outflows. We assume that the actual inflows from one component of the network (or from regions out of the network) to component $i\in\{1,...,n\}$ of the network are equal percentages of the corresponding attempted inflows, i.e.

$$\left(\begin{array}{c}\textit{actual flow of vehicles}\\ \textit{from component } j \textit{ to component } i\end{array}\right)=s_i(t)\left(\begin{array}{c}\textit{attempted flow of vehicles}\\ \textit{from component } j \textit{ to component } i\end{array}\right),\; i,j=1,...,n\,. \tag{3.8}$$

Other assumptions could be accommodated in the presented modeling framework if required. Combining equation (3.3) with equation (3.8) we get:

$$\left(\begin{array}{c}\textit{actual flow of vehicles}\\ \textit{from component } j \textit{ to component } i\end{array}\right)=s_i(t)p_{j,i}f_j(d,x_j)\,,\; i,j=1,...,n\,. \tag{3.9}$$

Moreover, we assume that the actual flow of vehicles from component $i\in\{1,...,n\}$ to regions out of the network is equal to the corresponding attempted flow of vehicles, which implies that there are no shock waves mounting from downstream and reaching the network exits. Thus, the actual outflow $q_i(t)\geq 0$ from component $i\in\{1,...,n\}$ of the network is given by:

$$q_i(t)=\left(Q_i+\sum_{j=1}^{n}s_j(t)p_{i,j}\right)f_i(d(t),x_i(t))\,,\; i=1,...,n\,,\; t\geq 0\,. \tag{3.10}$$

Combining equations (3.1), (3.2), (3.6), (3.7) and (3.10) we obtain the following discrete-time dynamical system:

$$x_i^+=x_i-\left(Q_i+\sum_{j=1}^{n}\frac{\min\left(a_j-x_j,v_j+\sum_{k=1}^{n}p_{k,j}f_k(d,x_k)\right)}{v_j+\sum_{k=1}^{n}p_{k,j}f_k(d,x_k)}p_{i,j}\right)f_i(d,x_i)+\min\left(a_i-x_i,v_i+\sum_{j=1}^{n}p_{j,i}f_j(d,x_j)\right),$$

$$\text{for } i=1,...,n \tag{3.11}$$

Define $S=[0,a_1]\times[0,a_2]\times\ldots\times[0,a_n]$. Since the functions $f_i:D\times[0,a_i]\to\Re_+$ satisfy $f_i(d,x_i)\leq x_i$ for all $(d,x_i)\in D\times[0,a_i]$, it follows that (3.11) is an uncertain dynamical system on $S=[0,a_1]\times[0,a_2]\times\ldots\times[0,a_n]$. We further assume that the uncertain dynamical system (3.11) has an equilibrium point. Specifically, we assume that:

**(H)** *The matrix* $P=\{p_{i,j}:i,j=1,...,n\}$ *satisfies* $\det(I-P')\neq 0$. *There exists a point* $x^*=(x_1^*,...,x_n^*)'\in S$ *that satisfies*

$$v_i+x_i^*+\sum_{j=1}^{n}p_{j,i}f_j(d,x_j^*)\leq a_i\text{, for all } d\in D \text{ and } i=1,...,n \tag{3.12}$$

$$f_i(d,x_i^*)=v_i+\sum_{j=1}^{n}p_{j,i}f_j(d,x_j^*)\text{, for all } d\in D \text{ and } i=1,...,n\,. \tag{3.13}$$

If we define the vector field $f(d,x)=(f_1(d,x_1),...,f_n(d,x_n))'$ and the vectors $v=(v_1,...,v_n)'\in\operatorname{int}(\Re_+^n)$, $a=(a_1,...,a_n)'\in\operatorname{int}(\Re_+^n)$, then conditions (3.12), (3.13) can be written in vector form:

$$(I-P')f(d,x^*)=v \text{ and } P'f(d,x^*)+v+x^*\leq a \tag{3.14}$$

and consequently we get $f(d,x^*) = f^* = (I-P')^{-1}v$ and $\left(I + P'(I-P')^{-1}\right)v + x^* \le a$. The equilibrium point $x^* = (x_1^*,...,x_n^*)' \in S$ is the so-called uncongested equilibrium point of the network. As mentioned in Section 1, the presented modeling framework is general enough to include several diverse kinds of network models and problems as special cases.

We are now in a position to state and prove the following theorem.

**Theorem 3.1:** *Consider system (3.11) under assumption (H). Assume that there exist constants* $0 \le b_i < c_i \le a_i$ *(* $i=1,...,n$ *) such that the set* $A = [b_1,c_1]\times[b_2,c_2]\times\ldots\times[b_n,c_n]$ *is a trapping region for system (3.11). Moreover, assume that there exists a constant* $L \ge 0$ *such that the following inequalities holds for all* $i=1,...,n$ *:*

$$\left|f_i(d,x_i) - f_i^*\right| \le L\left|x_i - x_i^*\right| \text{, for all } (d,x_i) \in D\times[0,a_i] \tag{3.15}$$

*Furthermore, assume that there exist constants* $\lambda_i, \mu_i \ge 0$, $\omega_i \in [x_i^*, a_i)$ *(* $i=1,...,n$ *) such that*

$$\left| x_i - x_i^* - \left( \left( Q_i + \sum_{j=1}^{n} \frac{\min\left( a_j - \theta_j\omega_j, v_j + p_{i,j}f_i(d,x_i) + \sum_{k\ne i} p_{k,j}f_k^* \right)}{v_j + p_{i,j}f_i(d,x_i) + \sum_{k\ne i} p_{k,j}f_k^*} p_{i,j} \right) f_i(d,x_i) + \min\left( a_i - x_i, v_i + \sum_{j=1}^{n} p_{j,i}f_j^* \right) \right) \right| \le \lambda_i\left|x_i - x_i^*\right|,$$

$$\text{for all } (d,x_i)\in D\times[b_i,c_i],\ \theta = (\theta_1,...,\theta_n)\in[0,1]^n \text{ and } i=1,...,n \tag{3.16}$$

$$\left|f_i^* - f_i(d,x_i)\right| \le \mu_i\left|x_i - x_i^*\right| \text{, for all } (d,x_i)\in D\times[b_i,c_i] \text{ and } i=1,...,n. \tag{3.17}$$

*Define* $F_i = \max\limits_{x_i\in[b_i,c_i],d\in D}\left(f_i(d,x_i)\right)$ *for* $i=1,...,n$ *and assume that*

$$f_j^* + p_{i,j}(F_i - f_i^*) \le a_j \text{ for all } i,j=1,...,n \tag{3.18}$$

*Define the matrix* $\Gamma = \{\gamma_{i,j} : i,j=1,...,n\}$ *by means of the formulae:*

$$\gamma_{i,i} := \lambda_i \text{, for } i=1,...,n \tag{3.19}$$

$$\gamma_{i,j} := \frac{F_i p_{i,j}\max(0, c_j - \omega_j)}{\left(f_j^* + p_{i,j}\left(F_i - f_i^*\right)\right)\left(c_j - x_j^*\right)} + \left( p_{j,i} + \sum_{k=1}^{n} \frac{F_i p_{i,k} p_{j,k}}{f_k^* + p_{i,k}\left(F_i - f_i^*\right)} \right)\mu_j \text{ for } i,j=1,...,n \text{ with } i\ne j. \tag{3.20}$$

*Suppose that the spectral radius* $\rho(\Gamma)$ *of the matrix* $\Gamma = \left\{\gamma_{i,j} \ge 0, i,j=1,...,n\right\} \in \Re_+^{n\times n}$ *is less than 1. Then the equilibrium point* $x^* = (x_1^*,...,x_n^*)' \in S$ *is RGES for system (3.11).*

**Remark 3.2:**

- **(a)** The constants $0 < c_i \le a_i$ ( $i=1,...,n$ ) are the ultimate upper bounds of the vehicle densities in each component of the network. In general, it is crucial to obtain accurate estimates of the constants $0 < c_i \le a_i$ ( $i=1,...,n$ ) (as small as possible) in order to show that congestion can only be a transient phenomenon and that there exist no congested equilibria for the network. The following section will also show the importance of the ultimate upper bounds of the vehicle densities in each component of the network.
- **(b)** Assumption (3.18) is not restrictive: since we are studying the stability properties of the uncongested equilibrium point, the equilibrium flow values $f_i^*$ for $i=1,...,n$ are far smaller than the capacities $a_i$ of the corresponding component, and condition (3.18) holds. In all practical examples that have been studied, condition (3.18) has been promptly satisfied.
- **(c)** Inequalities (3.16) are equivalent to the following inequalities for all $(d,x_i)\in D\times[b_i,c_i]$ and $i=1,...,n$ :

$$x_i^* - x_i + \left( Q_i + \sum_{j=1}^{n} \frac{\min\left( a_j , v_j + p_{i,j} f_i(d,x_i) + \sum_{k \neq i} p_{k,j} f_k^* \right)}{v_j + p_{i,j} f_i(d,x_i) + \sum_{k \neq i} p_{k,j} f_k^*} p_{i,j} \right) f_i(d,x_i) - \min\left( a_i - x_i , v_i + \sum_{j=1}^{n} p_{j,i} f_j^* \right) \le \lambda_i \left| x_i - x_i^* \right|$$

$$x_i - x_i^* - \left( Q_i + \sum_{j=1}^{n} \frac{\min\left( a_j - \omega_j , v_j + p_{i,j} f_i(d,x_i) + \sum_{k \neq i} p_{k,j} f_k^* \right)}{v_j + p_{i,j} f_i(d,x_i) + \sum_{k \neq i} p_{k,j} f_k^*} p_{i,j} \right) f_i(d,x_i) + \min\left( a_i - x_i , v_i + \sum_{j=1}^{n} p_{j,i} f_j^* \right) \le \lambda_i \left| x_i - x_i^* \right|$$

**Proof of Theorem 3.1:** We use Theorem 2.3 for the family of functions

$$V_i(x) := \left| x_i - x_i^* \right| \quad ( i = 1,...,n ) \tag{3.21}$$

and the dynamical system (3.11). Since the inequality

$$\frac{1}{\sqrt{n}} \left| x - x^* \right| \le \max_{i=1,...,n} \left( V_i(x) \right) \le \left| x - x^* \right|$$

holds for all $x \in A$ and since (3.15) implies the condition $\max\left\{ \left| F(d,x) - x^* \right| : d \in D \right\} \le \tilde{L} \left| x - x^* \right|$ for all $x \in S$, for certain constant $\tilde{L} \ge 0$, where $F(d,x) = (F_1(d,x),...,F_n(d,x))' \in \Re^n$ and

$$F_i(d,x) := x_i - \left( Q_i + \sum_{j=1}^{n} \frac{\min\left( a_j - x_j , v_j + \sum_{k=1}^{n} p_{k,j} f_k(d,x_k) \right)}{v_j + \sum_{k=1}^{n} p_{k,j} f_k(d,x_k)} p_{i,j} \right) f_i(d,x_i) + \min\left( a_i - x_i , v_i + \sum_{j=1}^{n} p_{j,i} f_j(d,x_j) \right)$$

it suffices to show that (2.3) holds for all $x \in A$, $i = 1,...,n$.

The rest part of proof is devoted to the proof of (2.3).

Indeed, using (3.21) we get for all $(d,x) \in D \times A$, $\theta = (\theta_1,...,\theta_n) \in [0,1]^n$ and $i = 1,...,n$:

$$\begin{aligned} V_i(F(d,x)) &= \left| x_i^+ - x_i^* \right| \\ &\le \left| x_i - x_i^* - \left( Q_i - \sum_{j=1}^{n} \frac{\min\left( a_j - \theta_j \omega_j , v_j + p_{i,j} f_i(d,x_i) + \sum_{k \neq i} p_{k,j} f_k^* \right)}{v_j + p_{i,j} f_i(d,x_i) + \sum_{k \neq i} p_{k,j} f_k^*} p_{i,j} \right) f_i(d,x_i) + \min\left( a_i - x_i , v_i + \sum_{j=1}^{n} p_{j,i} f_j^* \right) \right| \\ &\quad + f_i(d,x_i) \sum_{j=1}^{n} p_{i,j} w_{i,j} + \left| \min\left( a_i - x_i , v_i + \sum_{j=1}^{n} p_{j,i} f_j(d,x_j) \right) - \min\left( a_i - x_i , v_i + \sum_{j=1}^{n} p_{j,i} f_j^* \right) \right| \end{aligned} \tag{3.22}$$

where

$$w_{i,j} := \left| \frac{\min\left(a_j - x_j, v_j + \sum_{k=1}^{n} p_{k,j} f_k(d,x_k)\right)}{v_j + \sum_{k=1}^{n} p_{k,j} f_k(d,x_k)} - \frac{\min\left(a_j - \theta_j \omega_j, v_j + p_{i,j} f_i(d,x_i) + \sum_{k \neq i} p_{k,j} f_k^*\right)}{v_j + p_{i,j} f_i(d,x_i) + \sum_{k \neq i} p_{k,j} f_k^*} \right| . \tag{3.23}$$

Using (3.16), (3.17) and the fact that the inequality $\left|\min(a,x) - \min(a,y)\right| \le |x-y|$ holds for all $a,x,y \in \Re$, we obtain from (3.22) for all $(d,x) \in D \times A$, $\theta = (\theta_1,...,\theta_n) \in [0,1]^n$ and $i = 1,...,n$:

$$V_i(F(d,x)) \le \lambda_i \left| x_i - x_i^* \right| + f_i(d,x_i) \sum_{j=1}^{n} p_{i,j} w_{i,j} + \sum_{j=1}^{n} p_{j,i} \mu_j \left| x_j - x_j^* \right| . \tag{3.24}$$

We next show that for every $(d,x) \in D \times A$ and $i = 1,...,n$ we can select $\theta_j \in [0,1]$ in an appropriate way so that we can minimize the values of $w_{i,j}$ ($j = 1,...,n$).

Continuity of the mapping $[0,1] \ni \theta_j \to \min\left( \frac{a_j - \theta_j \omega_j}{v_j + p_{i,j} f_i(d,x_i) + \sum_{k \neq i} p_{k,j} f_k^*}, 1 \right)$ as well as the fact that the mapping $[0,1] \ni \theta_j \to \min\left( \frac{a_j - \theta_j \omega_j}{v_j + p_{i,j} f_i(d,x_i) + \sum_{k \neq i} p_{k,j} f_k^*}, 1 \right)$ is non-increasing, implies the existence of $\theta_j \in [0,1]$ with $w_{i,j} = 0$, provided that the following inequality holds:

$$\min\left( \frac{a_j - \omega_j}{v_j + p_{i,j} f_i(d,x_i) + \sum_{k \neq i} p_{k,j} f_k^*}, 1 \right) \le s_j \le \min\left( \frac{a_j}{v_j + p_{i,j} f_i(d,x_i) + \sum_{k \neq i} p_{k,j} f_k^*}, 1 \right) \tag{3.25}$$

where $s_j = \dfrac{\min\left(a_j - x_j, v_j + \sum_{k=1}^{n} p_{k,j} f_k(d,x_k)\right)}{v_j + \sum_{k=1}^{n} p_{k,j} f_k(d,x_k)}$. If inequality (3.25) does not hold, then we must have $\min\left( \frac{a_j - \omega_j}{v_j + p_{i,j} f_i(d,x_i) + \sum_{k \neq i} p_{k,j} f_k^*}, 1 \right) > s_j$. Indeed, this observation follows from the fact that $\min\left( \frac{a_j}{v_j + p_{i,j} f_i(d,x_i) + \sum_{k \neq i} p_{k,j} f_k^*}, 1 \right) = 1$ (a direct consequence of (3.13), (3.18)) and the fact that $s_j \le 1$.

Consequently, (3.23) implies that $w_{i,j}=\frac{\min\left(a_j-\omega_j, v_j+p_{i,j}f_i(d,x_i)+\sum_{k\neq i}p_{k,j}f_k^*\right)}{v_j+p_{i,j}f_i(d,x_i)+\sum_{k\neq i}p_{k,j}f_k^*}-\frac{\min\left(a_j-x_j, v_j+\sum_{k=1}^{n}p_{k,j}f_k(d,x_k)\right)}{v_j+\sum_{k=1}^{n}p_{k,j}f_k(d,x_k)}$ when (3.25) does not hold.

Moreover, since $s_j=\frac{\min\left(a_j-x_j, v_j+\sum_{k=1}^{n}p_{k,j}f_k(d,x_k)\right)}{v_j+\sum_{k=1}^{n}p_{k,j}f_k(d,x_k)}<1$, we get $a_j-x_j<v_j+\sum_{k=1}^{n}p_{k,j}f_k(d,x_k)$ and consequently,

$$w_{i,j}=\frac{\min\left(a_j-\omega_j, v_j+p_{i,j}f_i(d,x_i)+\sum_{k\neq i}p_{k,j}f_k^*\right)}{v_j+p_{i,j}f_i(d,x_i)+\sum_{k\neq i}p_{k,j}f_k^*}-\frac{a_j-x_j}{v_j+\sum_{k=1}^{n}p_{k,j}f_k(d,x_k)} \tag{3.26}$$

provided that (3.25) does not hold.

Hence, when (3.25) does not hold, we have from (3.26)

$$\begin{aligned} w_{i,j}=&\frac{\min\left(a_j-\omega_j, v_j+p_{i,j}f_i(d,x_i)+\sum_{k\neq i}p_{k,j}f_k^*\right)-(a_j-x_j)}{v_j+p_{i,j}f_i(d,x_i)+\sum_{k\neq i}p_{k,j}f_k^*} \\ &+\frac{a_j-x_j}{v_j+p_{i,j}f_i(d,x_i)+\sum_{k\neq i}p_{k,j}f_k^*}-\frac{a_j-x_j}{v_j+\sum_{k=1}^{n}p_{k,j}f_k(d,x_k)} \end{aligned}$$

(manipulation)

$$\begin{aligned} w_{i,j}\leq&\frac{\min\left(a_j-\omega_j, v_j+p_{i,j}f_i(d,x_i)+\sum_{k\neq i}p_{k,j}f_k^*\right)-(a_j-x_j)}{v_j+p_{i,j}f_i(d,x_i)+\sum_{k\neq i}p_{k,j}f_k^*} \\ &+\frac{a_j-x_j}{\left(v_j+p_{i,j}f_i(d,x_i)+\sum_{k\neq i}p_{k,j}f_k^*\right)\left(v_j+\sum_{k=1}^{n}p_{k,j}f_k(d,x_k)\right)}\sum_{k\neq i}p_{k,j}\mu_k\left|x_k-x_k^*\right| \end{aligned}$$

(Use of (3.17))

$$\begin{aligned} &w_{i,j}\left(v_j+p_{i,j}f_i(d,x_i)+\sum_{k\neq i}p_{k,j}f_k^*+\sum_{k\neq i}p_{k,j}\mu_k\left|x_k-x_k^*\right|\right)\leq\min\left(a_j-\omega_j, v_j+p_{i,j}f_i(d,x_i)+\sum_{k\neq i}p_{k,j}f_k^*\right)-(a_j-x_j) \\ &+\frac{\min\left(a_j-\omega_j, v_j+p_{i,j}f_i(d,x_i)+\sum_{k\neq i}p_{k,j}f_k^*\right)}{v_j+p_{i,j}f_i(d,x_i)+\sum_{k\neq i}p_{k,j}f_k^*}\sum_{k\neq i}p_{k,j}\mu_k\left|x_k-x_k^*\right| \end{aligned}$$

(Use of (3.25))

$$w_{i,j} \le \frac{\max(0, x_j - \omega_j)}{v_j + p_{i,j} f_i(d,x_i) + \sum_{k \ne i} p_{k,j} f_k^* + \sum_{k \ne i} p_{k,j} \mu_k \left| x_k - x_k^* \right|}$$

$$+ \frac{\min\left( a_j - \omega_j, v_j + p_{i,j} f_i(d,x_i) + \sum_{k \ne i} p_{k,j} f_k^* \right) \sum_{k \ne i} p_{k,j} \mu_k \left| x_k - x_k^* \right|}{\left( v_j + p_{i,j} f_i(d,x_i) + \sum_{k \ne i} p_{k,j} f_k^* \right)\left( v_j + p_{i,j} f_i(d,x_i) + \sum_{k \ne i} p_{k,j} f_k^* + \sum_{k \ne i} p_{k,j} \mu_k \left| x_k - x_k^* \right| \right)}$$

(Use of $\min\left( a_j - \omega_j, v_j + p_{i,j} f_i(d,x_i) + \sum_{k \ne i} p_{k,j} f_k^* \right) \le a_j - \omega_j$ and $x_j - \omega_j \le \max(0, x_j - \omega_j)$)

$$w_{i,j} \le \frac{\max(0, x_j - \omega_j) + \sum_{k \ne i} p_{k,j} \mu_k \left| x_k - x_k^* \right|}{f_j^* + p_{i,j}\left( f_i(d,x_i) - f_i^* \right) + \sum_{k \ne i} p_{k,j} \mu_k \left| x_k - x_k^* \right|}$$

(Use of (3.13) and $\frac{\min\left( a_j - \omega_j, v_j + p_{i,j} f_i(d,x_i) + \sum_{k \ne i} p_{k,j} f_k^* \right)}{v_j + p_{i,j} f_i(d,x_i) + \sum_{k \ne i} p_{k,j} f_k^*} \le 1$)

The above inequality holds when (3.25) holds as well. Using the above inequality in conjunction with (3.24), we obtain for all $(d,x) \in D \times A$ and $i = 1,...,n$:

$$V_i(F(d,x)) \le \lambda_i \left| x_i - x_i^* \right| + f_i(d,x_i) \sum_{j=1}^{n} p_{i,j} \frac{\max(0, x_j - \omega_j) + \sum_{k \ne i} p_{k,j} \mu_k \left| x_k - x_k^* \right|}{f_j^* + \sum_{k \ne i} p_{k,j} \mu_k \left| x_k - x_k^* \right| + p_{i,j}\left( f_i(d,x_i) - f_i^* \right)} + \sum_{j=1}^{n} p_{j,i} \mu_j \left| x_j - x_j^* \right|. \qquad (3.27)$$

Using the facts that $x_j \in [b_j, c_j]$ and $\omega_j \ge x_j^*$, we obtain $\max(0, x_j - \omega_j) \le \frac{\max(0, c_j - \omega_j)}{c_j - x_j^*} \left| x_j - x_j^* \right|$. Therefore, we obtain from (3.27) for all $(d,x) \in D \times A$ and $i = 1,...,n$:

$$\begin{aligned} V_i(F(d,x)) \le \lambda_i \left| x_i - x_i^* \right| &+ \sum_{j=1}^{n} \frac{p_{i,j} f_i(d,x_i) \max(0, c_j - \omega_j)}{\left( f_j^* + p_{i,j}\left( f_i(d,x_i) - f_i^* \right) \right)\left( c_j - x_j^* \right)} \left| x_j - x_j^* \right| + \sum_{j=1}^{n} p_{j,i} \mu_j \left| x_j - x_j^* \right| \\ &+ \sum_{j=1}^{n} \frac{p_{i,j} f_i(d,x_i)}{f_j^* + p_{i,j}\left( f_i(d,x_i) - f_i^* \right)} \sum_{k \ne i} p_{k,j} \mu_k \left| x_k - x_k^* \right| \end{aligned} \qquad (3.28)$$

Finally, using definitions (3.21), the fact that $\frac{p_{i,j} f_i(d,x_i)}{f_j^* + p_{i,j}\left( f_i(d,x_i) - f_i^* \right)} \le \frac{p_{i,j} F_i}{f_j^* + p_{i,j}\left( F_i - f_i^* \right)}$ for all $(d,x) \in D \times A$, where $F_i = \max_{x_i \in [b_i, c_i], d \in D} \left( f_i(d,x_i) \right)$ for $i = 1,...,n$ and the fact that $p_{i,i} = 0$, we obtain (2.3).

The proof is complete. $\lhd$

**Example 3.3:** Consider the traffic network shown in Figure 1, for which the matrix $P=\{p_{i,j}:i,j=1,...,5\}$ is:

$$P=\begin{bmatrix} 0 & p & 0 & 0 & 0 \\ 0 & 0 & p & 0 & \tilde{p} \\ p & 0 & 0 & \tilde{p} & 0 \\ 0 & 0 & 0 & 0 & 0 \\ 0 & 0 & 0 & 0 & 0 \end{bmatrix} \tag{3.29}$$

where $p,\tilde{p}>0$ are constants with $p+\tilde{p}\le 1$. The external inflows and the capacities of the network are:

$$v_1=v_2=v_3=v>0\,,\ v_4=v_5=\tilde{v}>0\,,\ a_i=a>0\ (i=1,...,5) \tag{3.30}$$

where $v,\tilde{v},a>0$ are constants. Finally, we assume that all functions $f_i:D\times[0,a]\to\Re_+$ ($i=1,...,5$) are independent of $d\in D$ and satisfy:

$$f_i(x)=f(x):=\begin{cases} rx\ for\ x\in[0,\delta] \\ r\delta-q(x-\delta)\ for\ x\in(\delta,a] \end{cases}\ (i=1,...,5) \tag{3.31}$$

where $\delta\in(0,a)$, $r\in(0,1]$, $q\in\left[0,\frac{\delta}{a-\delta}r\right]$ are constants. Note that the lower part of the right-hand side of (3.31) allows for the modeling of capacity drop at the outflow of congestion according to [15]. The network has the (uncongested) equilibrium point

$$x^*=\left(\frac{v}{r(1-p)},\frac{v}{r(1-p)},\frac{v}{r(1-p)},\frac{\tilde{v}(1-p)+\tilde{p}v}{r(1-p)},\frac{\tilde{v}(1-p)+\tilde{p}v}{r(1-p)}\right) \tag{3.32}$$

which satisfies assumption (H) provided that $\frac{v}{1-p}\le r\delta$, $\frac{v}{1-p}\left(1+\frac{1}{r}\right)\le a$, $\frac{\tilde{v}(1-p)+\tilde{p}v}{1-p}\le r\delta$ and $\left(1+\frac{1}{r}\right)\frac{\tilde{v}(1-p)+\tilde{p}v}{1-p}\le a$. We next apply Theorem 3.1 under the assumption

$$v+pr\delta\le a\ \text{and}\ \tilde{v}+\tilde{p}r\delta\le a \tag{3.33}$$

with $A=[0,a]^5$. Assumption (3.33) is the expression of assumption (3.18) for the given network. The matrix $\Gamma=\{\gamma_{i,j}\ge 0,i,j=1,...,5\}\in\Re_+^{5\times 5}$ is equal to:

$$\Gamma=\begin{bmatrix} \lambda_1 & (a-\omega_2)\varphi & p\mu & 0 & 0 \\ p\mu & \lambda_2 & (a-\omega_3)\varphi & 0 & (a-\omega_5)\zeta \\ (a-\omega_1)\varphi & p\mu & \lambda_3 & (a-\omega_4)\zeta & 0 \\ 0 & 0 & \tilde{p}\mu & \lambda_4 & 0 \\ 0 & \tilde{p}\mu & 0 & 0 & \lambda_5 \end{bmatrix} \tag{3.34}$$

where

$$\varphi=\frac{pr^2\delta(1-p)}{(pr\delta+v)(ar(1-p)-v)},\ \zeta=\frac{r^2\delta\tilde{p}(1-p)}{(\tilde{v}+\tilde{p}r\delta)(ar(1-p)-\tilde{v}(1-p)-\tilde{p}v)} \tag{3.35}$$

$$\mu:=\sup\left\{r\frac{|v-(1-p)f(s)|}{|rs(1-p)-v|}:s\in[0,a],s\neq\frac{v}{r(1-p)}\right\} \tag{3.36}$$

$$\lambda_4 = \lambda_5 = \sup\left\{ \frac{\left| s - \frac{\tilde{v}(1-p)+\tilde{p}v}{r(1-p)} - f(s) + \min\left( a - s, \frac{\tilde{v}(1-p)+\tilde{p}v}{1-p} \right) \right|}{\left| s - \frac{\tilde{v}(1-p)+\tilde{p}v}{r(1-p)} \right|} : s \in [0,a], s \neq \frac{\tilde{v}(1-p)+\tilde{p}v}{r(1-p)} \right\} \tag{3.37}$$

$$\lambda_i = \max(u_i, w_i)\,,\ i = 1,2,3 \tag{3.38}$$

where

$$u_1 := \sup\left\{ \frac{s - \frac{v}{r(1-p)} - \left( 1 - p + \frac{\min(a-\omega_2, v+pf(s))}{v+pf(s)} p \right) f(s) + \min\left( a-s, \frac{v}{1-p} \right)}{\left| s - \frac{v}{r(1-p)} \right|} : s \in [0,a], s \neq \frac{v}{r(1-p)} \right\}$$

$$w_1 = \sup\left\{ \frac{\frac{v}{r(1-p)} - s + \left( 1 - p + \frac{\min(a, v+pf(s))}{v+pf(s)} p \right) f(s) - \min\left( a-s, \frac{v}{1-p} \right)}{\left| s - \frac{v}{r(1-p)} \right|} : s \in [0,a], s \neq \frac{v}{r(1-p)} \right\}$$

$$u_2 = \sup\left\{ \frac{s - \frac{v}{r(1-p)} - \left( 1 - p - \tilde{p} + \frac{\min(a-\omega_3, v+pf(s))}{v+pf(s)} p + \frac{\min(a-\omega_5, \tilde{v}+\tilde{p}f(s))}{\tilde{v}+\tilde{p}f(s)} \tilde{p} \right) f(s) + \min\left( a-s, \frac{v}{1-p} \right)}{\left| s - \frac{v}{r(1-p)} \right|} : s \in [0,a], s \neq \frac{v}{r(1-p)} \right\}$$

$$w_2 = w_3 = \sup\left\{ \frac{\frac{v}{r(1-p)} - s + \left( 1 - p - \tilde{p} + \frac{\min(a, v+pf(s))}{v+pf(s)} p + \frac{\min(a, \tilde{v}+\tilde{p}f(s))}{\tilde{v}+\tilde{p}f(s)} \tilde{p} \right) f(s) - \min\left( a-s, \frac{v}{1-p} \right)}{\left| s - \frac{v}{r(1-p)} \right|} : s \in [0,a], s \neq \frac{v}{r(1-p)} \right\}$$

$$u_3 = \sup\left\{ \frac{s - \frac{v}{r(1-p)} - \left( 1 - p - \tilde{p} + \frac{\min(a-\omega_1, v+pf(s))}{v+pf(s)} p + \frac{\min(a-\omega_4, \tilde{v}+\tilde{p}f(s))}{\tilde{v}+\tilde{p}f(s)} \tilde{p} \right) f(s) + \min\left( a-s, \frac{v}{1-p} \right)}{\left| s - \frac{v}{r(1-p)} \right|} : s \in [0,a], s \neq \frac{v}{r(1-p)} \right\}$$

and $\omega_i \in \left[ \frac{v}{r(1-p)}, a \right)$ $(i = 1,2,3)$, $\omega_i \in \left[ \frac{\tilde{v}(1-p)+\tilde{p}v}{r(1-p)}, a \right)$ $(i = 4,5)$ are constants. All quantities involved in (3.35)-(3.38) can be computed numerically once the constants $p, \tilde{p} > 0$, $v, \tilde{v}, a > 0$, $\delta \in (0,a)$, $r \in (0,1]$, $q \in \left[ 0, \frac{\delta}{a-\delta} r \right]$ and $\omega_i \in [0,a)$ $(i = 1,...,5)$ are known. More specifically, the constants $\omega_i \in \left[ \frac{v}{r(1-p)}, a \right)$ $(i = 1,2,3)$, $\omega_i \in \left[ \frac{\tilde{v}(1-p)+\tilde{p}v}{r(1-p)}, a \right)$ $(i = 4,5)$ can be selected in an appropriate way so that the spectral radius of $\Gamma$ is minimized.

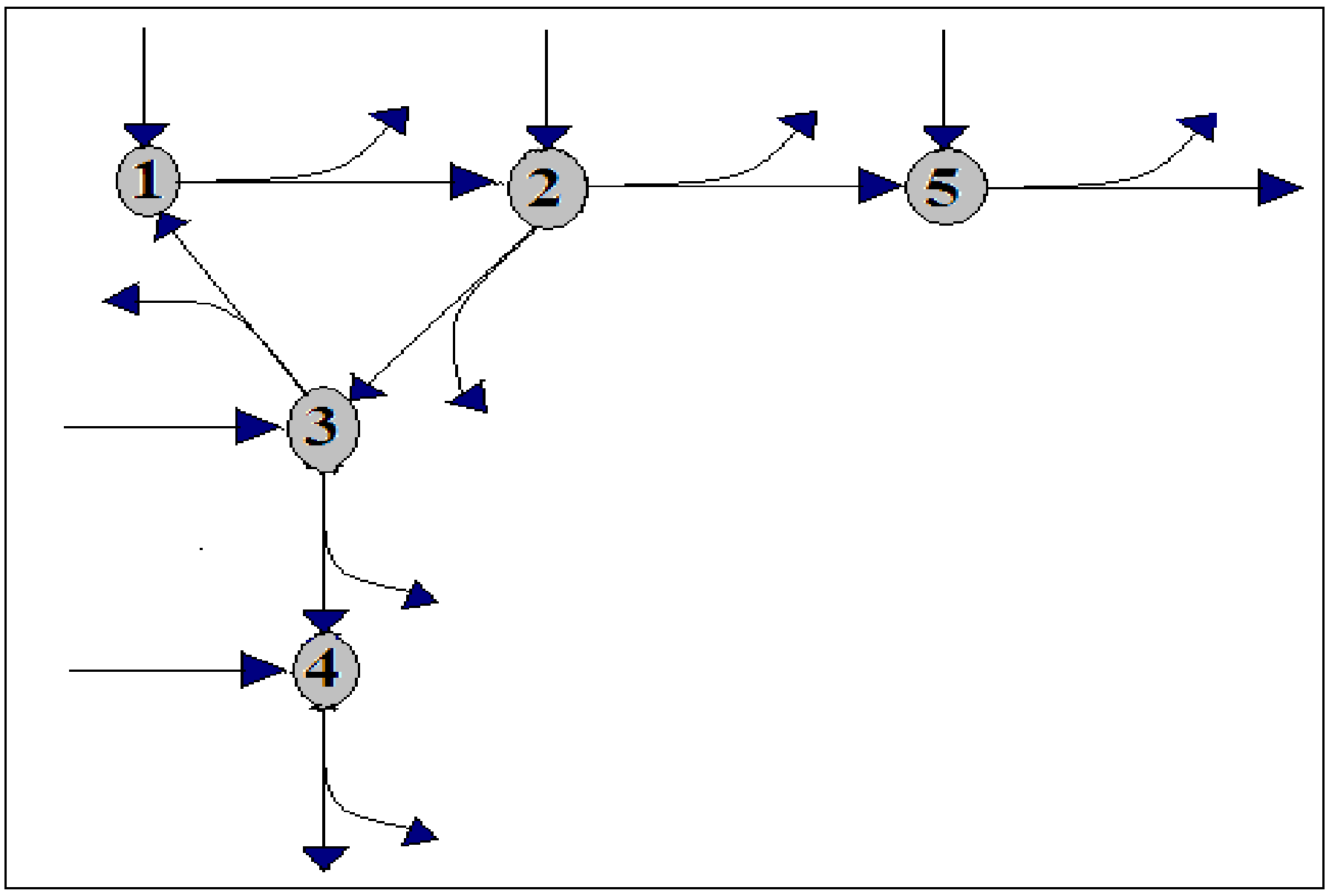

**Fig. 1:** The traffic network of Example 3.3

For $a=10$, $v=0.4$, $\tilde{v}=0.4$, $\delta=5$, $p=0.2$, $\tilde{p}=0.1$, $r=0.55$, $q=0.1$, the selection $\omega_1=9.14$, $\omega_2=8.532$, $\omega_3=9.559$, $\omega_4=9.3697$, $\omega_5=9.329$ gives us the matrix $\Gamma=\{\gamma_{i,j}\geq 0, i,j=1,...,5\}\in\Re_+^{5\times 5}$:

$$\Gamma=\begin{bmatrix} 0.7905 & 0.0281 & 0.11 & 0 & 0 \\ 0.11 & 0.8166 & 0.0281 & 0 & 0.0298 \\ 0.0548 & 0.11 & 0.7905 & 0.028 & 0 \\ 0 & 0 & 0.055 & 0.7869 & 0 \\ 0 & 0.055 & 0 & 0 & 0.7869 \end{bmatrix}$$

Since $\max\limits_{i=1,...,5}\left(\sum_{j=1}^{5}\gamma_{i,j}\right)=0.9845<1$, we can conclude that $\rho(\Gamma)<1$ and consequently, Theorem 3.1 implies that the (uncongested) equilibrium point is GES. $\lhd$

## 4. Global Exponential Stability Results for Freeways

A freeway divided in $n\geq 3$ sections or cells, each with one on-ramp and one off-ramp, is a traffic network of the form (3.11) with $p_{i,j}=0$ for all $i,j=1,...,n$ with $j\neq i+1$. Defining $p_{i,i+1}=p_i$ for $i=1,...,n-1$, we obtain from (3.5) and (3.11) the model:

$$x_1^+=x_1-\left(1-p_1+\frac{\min(a_2-x_2,v_2+p_1f_1(d,x_1))}{v_2+p_1f_1(d,x_1)}p_1\right)f_1(d,x_1)+\min(a_1-x_1,v_1) \tag{4.1}$$

$$x_i^+=x_i-\left(1-p_i+\frac{\min(a_{i+1}-x_{i+1},v_{i+1}+p_if_i(d,x_i))}{v_{i+1}+p_if_i(d,x_i)}p_i\right)f_i(d,x_i)+\min(a_i-x_i,v_i+p_{i-1}f_{i-1}(d,x_{i-1})),$$
$$\text{for } i=2,...,n-1 \tag{4.2}$$

$$x_n^+=x_n-f_n(d,x_n)+\min(a_n-x_n,v_n+p_{n-1}f_{n-1}(d,x_{n-1})). \tag{4.3}$$

If we further suppose that $v_i = 0$ for $i = 2,...,n$, $p_i = 1$, $f_i(d, x_i) = f_i(x_i)$ for $i = 1,...,n$, $v_1 = v > 0$, the model (4.1), (4.2), (4.3) becomes the following disturbance-free model:

$$x_1^+ = x_1 - \min(a_2 - x_2, f_1(x_1)) + \min(a_1 - x_1, v) \tag{4.4}$$

$$x_i^+ = x_i - \min(a_{i+1} - x_{i+1}, f_i(x_i)) + \min(a_i - x_i, f_{i-1}(x_{i-1})), \text{ for } i = 2,...,n-1 \tag{4.5}$$

$$x_n^+ = x_n - f_n(x_n) + \min(a_n - x_n, f_{n-1}(x_{n-1})). \tag{4.6}$$

Again $f_i : [0, a_i] \to \Re_+$ ($i = 1,...,n$) are continuous functions with $f_i(s) \le s$ for all $s \in [0, a_i]$. We suppose that there exists a vector $x^* = (x_1^*,...,x_n^*) \in [0, a_1] \times \ldots [0, a_n]$ with $f_i(x_i^*) = v$ and $x_i^* + v < a_i$ ($i = 1,...,n$). It follows that assumption (H) holds for the equilibrium point $x^* \in \Re^n$. The following corollary is a direct consequence of Theorem 3.1 (although Theorem 3.1 was applied to the model (3.11) which required $v_i > 0$ for $i = 1,...,n$, all arguments in the proof of Theorem 3.1 can be repeated).

**Corollary 4.1:** *Consider system (4.4), (4.5), (4.6) with $n \ge 3$. Assume that there exist constants $0 \le b_i < c_i \le a$ ($i = 1,...,n$) such that the set $A = [b_1, c_1] \times [b_2, c_2] \times \ldots \times [b_n, c_n]$ is a trapping region for system (4.4), (4.5), (4.6). Moreover, assume that there exists a constant $L \ge 0$ such that the following inequalities hold for $i = 1,...,n$:*

$$|f_i(x) - v| \le L|x - x_i^*|, \text{ for all } x \in [0, a_i]. \tag{4.7}$$

*Furthermore, assume that there exist constants $\lambda_i \ge 0$ ($i = 1,...,n$), $\mu_i \ge 0$ ($i = 1,...,n-1$), $\omega_i \in [x_i^*, a_i)$ ($i = 2,...,n$) such that*

$$\left|s - x_i^* - \min(a_{i+1} - \omega_{i+1}, f_i(s)) + \min(a_i - s, v)\right| \le \lambda_i |s - x_i^*|, \text{ for } s \in [b_i, c_i],\ i = 1,...,n-1 \tag{4.8}$$

$$\left|s - x_i^* - \min(a_{i+1}, f_i(s)) + \min(a_i - s, v)\right| \le \lambda_i |s - x_i^*|, \text{ for } s \in [b_i, c_i],\ i = 1,...,n-1 \tag{4.9}$$

$$\left|s - x_n^* - f_n(s) + \min(a_n - s, v)\right| \le \lambda_n |s - x_n^*|, \text{ for all } s \in [b_n, c_n] \tag{4.10}$$

$$|v - f_i(s)| \le \mu_i |s - x_i^*|, \text{ for all } s \in [b_i, c_i],\ i = 1,...,n-1. \tag{4.11}$$

*Assume that $\max_{s \in [b_i, c_i]} (f_i(s)) \le a_{i+1}$ for all $i = 1,...,n-1$. Define the tridiagonal matrix $\Gamma = \{\gamma_{i,j} : i, j = 1,...,n\}$ by means of the formulae:*

$$\gamma_{i,i} := \lambda_i, \text{ for } i = 1,...,n \tag{4.12}$$

$$\gamma_{i,i+1} := \frac{\max(0, c_{i+1} - \omega_{i+1})}{c_{i+1} - x_i^*} \text{ for } i = 1,...,n-1 \tag{4.13}$$

$$\gamma_{i,i-1} := \mu_{i-1}, \text{ for } i = 2,...,n. \tag{4.14}$$

*Suppose that the spectral radius $\rho(\Gamma)$ of the matrix $\Gamma = \{\gamma_{i,j} \ge 0, i, j = 1,...,n\} \in \Re_+^{n \times n}$ is less than 1. Then the equilibrium point $x^* = (x_1^*,...,x_n^*) \in [0, a_1] \times \ldots [0, a_n]$ is GES for the disturbance-free system (4.4), (4.5), (4.6).*

Corollary 4.1 shows more emphatically than the general Theorem 3.1 that the determination of a trapping region is crucial for the stability properties of the disturbance-free system (4.4), (4.5), (4.6). Indeed, if $\omega_i \ge c_i$ for $i = 2,...,n$ then the spectral radius $\rho(\Gamma)$ of the matrix $\Gamma = \{\gamma_{i,j} \ge 0, i, j = 1,...,n\} \in \Re_+^{n \times n}$ is simply $\rho(\Gamma) = \max_{i=1,...,n} (\lambda_i)$ (because in this case the matrix $\Gamma$ is lower triangular). Of course, the crudest trapping region that can be used is the set $A = [0, a_1] \times \ldots [0, a_n]$. However, we can generate "smaller" trapping regions by means of the following propositions.

**Proposition 4.2:** *Suppose that that there exist constants* $0 \le b_i < c_i \le a_i$ $(i = 1,...,n)$ *such that the set* $A = [b_1, c_1] \times [b_2, c_2] \times \ldots \times [b_n, c_n]$ *is a trapping region for system (4.4), (4.5), (4.6) with* $n \ge 3$. *Let* $i \in \{1,...,n\}$, $c \in [0, c_i]$ *be a constant such that one of the following implications hold:*

$$\text{If } i = 1 \text{ and } c \ge x_1^* \text{ then } \min_{c \le s \le c_1} \left( \min(a_2 - c_2, f_1(s)) - \min(a_1 - s, v) \right) > 0$$
$$\text{and } \max_{b_1 \le s \le c} \left( s - \min(a_2 - c_2, f_1(s)) + \min(a_1 - s, v) \right) \le c. \tag{4.15}$$

$$\text{If } i \in \{2,...,n-1\} \text{ and } c \ge x_i^* \text{ then } \min_{c \le s \le c_i} \left( \min(a_{i+1} - c_{i+1}, f_i(s)) - \min(a_i - s, F_{i-1}) \right) > 0 \text{ and}$$
$$\max_{b_i \le s \le c} \left( s - \min(a_{i+1} - c_{i+1}, f_i(s)) + \min(a_i - s, F_{i-1}) \right) \le c, \text{ where } F_{i-1} := \max_{s \in [b_{i-1}, c_{i-1}]} \left( f_{i-1}(s) \right). \tag{4.16}$$

$$\text{If } i = n \text{ and } c \ge x_n^* \text{ then } \min_{c \le s \le c_n} \left( f_n(s) - \min(a_n - s, F_{n-1}) \right) > 0 \text{ and}$$
$$\max_{b_n \le s \le c} \left( s - f_n(s) + \min(a_n - s, F_{n-1}) \right) \le c, \text{ where } F_{n-1} := \max_{s \in [b_{n-1}, c_{n-1}]} \left( f_{n-1}(s) \right). \tag{4.17}$$

*Then the set* $B \subseteq A$ *which results from the replacement of* $c_i$ *by* $c$ *in the formula* $[b_1, c_1] \times [b_2, c_2] \times \ldots \times [b_n, c_n]$ *is a trapping region for system (4.4), (4.5), (4.6).*

**Proof:** We consider the case $i = 1$ (all other cases are similar). We only consider the case $c < c_1$ (since the case $c = c_1$, is trivial). Notice that since $A = [b_1, c_1] \times [b_2, c_2] \times \ldots \times [b_n, c_n]$ is a trapping region for system (4.4), (4.5), (4.6), there exists $m \ge 0$ such that for every $x_0 \in S = [0, a_1] \times \ldots [0, a_n]$ the solution $x(t)$ of (4.4), (4.5), (4.6) with initial condition $x(0) = x_0$ satisfies $x(t) \in A$ for all $t \ge m$. Consequently, we have from (4.4) for all $t \ge m$:

$$x_1(t+1) \le x_1(t) - \min(a_2 - c_2, f_1(x_1(t))) + \min(a_1 - x_1(t), v). \tag{4.18}$$

It follows from (4.15) that, if $x_1(t) \le c$ for certain $t \ge m$ then $x_1(t+1) \le c$. Consequently, the following property holds:

**(P):** If there exists $T \ge m$ with $x_1(T) \le c$ then it holds that $x_1(t) \le c$ for all $t \ge T$.

Let $\delta := \min_{c \le s \le c_1} \left( \min(a_2 - c_2, f_1(s)) - \min(a_1 - s, v) \right) > 0$. We claim that the solution $x(t)$ of (4.4), (4.5), (4.6) with arbitrary initial condition $x(0) = x_0 \in S = [0, a_1] \times \ldots [0, a_n]$ satisfies $x_1(t) \le c$ for all $t \ge m + \left\lceil \frac{c_1 - c}{\delta} \right\rceil + 1$. The proof is made by contradiction. Suppose that there exists $x_0 \in S = [0, a_1] \times \ldots [0, a_n]$ and $t \ge m + \left\lceil \frac{c_1 - c}{\delta} \right\rceil + 1$ such that $x_1(t) > c$. Notice that property (P) guarantees that $x_1(j) > c$ for all $j = m,...,t$. It follows from (4.18) and definition $\delta := \min_{c \le s \le c_1} \left( \min(a_2 - c_2, f_1(s)) - \min(a_1 - s, v) \right) > 0$ that the following inequality holds for all $j = m,...,t$:

$$x_1(j+1) \le x_1(j) - \delta. \tag{4.19}$$

Inequality (4.19) implies that $x_1(t) \le x_1(m) - (t-m)\delta$. The previous inequality in conjunction with $x_1(t) > c$ and the fact that $x_1(m) \le c_1$ implies $(t-m)\delta < c_1 - c$ which contradicts the fact that $t \ge m + \left\lceil \frac{c_1 - c}{\delta} \right\rceil + 1$. The proof is complete. ◁

**Proposition 4.3:** *Suppose that there exist constants* $0 \le b_i < c_i \le a_i$ $(i = 1,...,n)$ *such that the set* $A = [b_1, c_1] \times [b_2, c_2] \times \ldots \times [b_n, c_n]$ *is a trapping region for system (4.4), (4.5), (4.6) with* $n \ge 3$. *Let* $i \in \{1,...,n\}$, $b \in [b_i, a_i]$ *be a constant such that one of the following implications hold:*

If $i=1$ and $b \le x_1^*$ then $v > \max_{b_1 \le s \le b}\left(f_1(s)\right)$ and $\min_{b \le s \le c_1}\left(s - \min\left(a_2 - b_2, f_1(s)\right) + \min\left(a_1 - s, v\right)\right) \ge b$. (4.20)

If $i \in \{2,...,n-1\}$ and $b \le x_i^*$ then $\max_{b_i \le s \le b}\left(\min\left(a_{i+1} - b_{i+1}, f_i(s)\right) - \min\left(a_i - s, \tilde{f}_{i-1}\right)\right) > 0$

and $\min_{b \le s \le c_i}\left(s - \min\left(a_{i+1} - b_{i+1}, f_i(s)\right) + \min\left(a_i - s, \tilde{f}_{i-1}\right)\right) \ge b$, where $\tilde{f}_{i-1} := \min_{s \in [b_{i-1}, c_{i-1}]}\left(f_{i-1}(s)\right)$. (4.21)

If $i=n$ and $b \le x_n^*$ then $\max_{b_n \le s \le b}\left(f_n(s) - \min\left(a_n - s, \tilde{f}_{n-1}\right)\right) > 0$

and $\min_{b \le s \le c_i}\left(s - f_n(s) + \min\left(a_n - s, \tilde{f}_{n-1}\right)\right) \ge b$, where $\tilde{f}_{n-1} := \min_{s \in [b_{n-1}, c_{n-1}]}\left(f_{n-1}(s)\right)$. (4.22)

*Then the set $B \subseteq A$ which results from the replacement of $b_i$ by $b$ in the formula $[b_1, c_1] \times [b_2, c_2] \times \ldots \times [b_n, c_n]$ is a trapping region for system (4.4), (4.5), (4.6).*

**Proof:** We consider the case $i=1$ (all other cases are similar). We only consider the case $b > b_1$ (since the case $b = b_1$, is trivial). Notice that since $A = [b_1, c_1] \times [b_2, c_2] \times \ldots \times [b_n, c_n]$ is a trapping region for system (4.4), (4.5), (4.6), there exists $m \ge 0$ such that for every $x_0 \in S = [0, a_1] \times \ldots [0, a_n]$ the solution $x(t)$ of (4.4), (4.5), (4.6) with initial condition $x(0) = x_0$ satisfies $x(t) \in A$ for all $t \ge m$. Consequently, we have from (4.4) for all $t \ge m$:

$$x_1(t+1) \ge x_1(t) - \min\left(a_2 - b_2, f_1(x_1(t))\right) + \min\left(a_1 - x_1(t), v\right) \tag{4.23}$$

It follows from (4.20) that, if $x_1(t) \ge b$ for certain $t \ge m$ then $x_1(t+1) \ge b$. Consequently, the following property holds:

**(P):** If there exists $T \ge m$ with $x_1(T) \ge b$ then it holds that $x_1(t) \ge b$ for all $t \ge T$.

Let $\delta := v - \max_{0 \le s \le b_1}\left(f_1(s)\right) > 0$. We claim that the solution $x(t)$ of (4.4), (4.5), (4.6) with arbitrary initial condition $x(0) = x_0 \in S = [0, a_1] \times \ldots [0, a_n]$ satisfies $x_1(t) \ge b$ for all $t \ge m + \left\lceil \frac{b - b_1}{\delta} \right\rceil + 1$. The proof is made by contradiction. Suppose that there exists $x_0 \in S = [0, a_1] \times \ldots [0, a_n]$ and $t \ge m + \left\lceil \frac{b - b_1}{\delta} \right\rceil + 1$ such that $x_1(t) < b$. Notice that property (P) guarantees that $x(j) < b$ for all $j = m,...,t$. It follows from (4.23) and the fact that $b \le x_1^* < a_1 - v$ that the following inequality holds for all $j = m,...,t$:

$$x_1(t+1) \ge x_1(t) - \min\left(a_2 - b_2, f_1(x_1(t))\right) + v \tag{4.24}$$

Since every trapping region must contain the equilibrium point of (4.4), (4.5), (4.6), it follows that $b_2 \le x_2^* < a_2 - v$ and since $\delta := v - \max_{b_1 \le s \le b}\left(f_1(s)\right) > 0$, we get $f_1(s) < v < a_2 - b_2$ for all $s \in [b_1, b]$. It follows from definition $\delta := v - \max_{b_1 \le s \le b}\left(f_1(s)\right) > 0$ that the following inequality holds for all $j = m,...,t$:

$$x_1(j+1) \ge x_1(j) + \delta \tag{4.25}$$

Inequality (4.25) implies that $x_1(t) \ge x_1(m) + (t-m)\delta$. The previous inequality in conjunction with $x_1(t) < b$ and the fact that $x_1(m) \ge b_1$ implies $(t-m)\delta < b - b_1$ which contradicts the fact that $t \ge m + \left\lceil \frac{b - b_1}{\delta} \right\rceil + 1$. The proof is complete. $\lhd$

Using Proposition 4.2 and Corollary 4.1, we can construct an algorithm that provides easily checkable sufficient conditions for the GES of the equilibrium point $x^* = (x_1^*,...,x_n^*) \in [0, a_1] \times \ldots [0, a_n]$ of the disturbance-free system (4.4), (4.5), (4.6).

**Corollary 4.4:** *Consider system (4.4), (4.5), (4.6) with* $n \geq 3$. *Suppose that* $0 < f_i(s) < a_{i+1}$ *for all* $s \in (0, a_i]$ *and* $i = 1,...,n-1$. *Perform the following algorithm:*

Step 1: Find $k_n \in [x_n^*, a_n)$ such that $\min_{k_n \leq s \leq a_n} \left( f_n(s) - \min(a_n - s, F_{n-1}) \right) > 0$ and $\max_{0 \leq s \leq k_n} \left( s - f_n(s) + \min(a_n - s, F_{n-1}) \right) \leq k_n$, where $F_{n-1} := \max_{s \in [0, a_{n-1}]} \left( f_{n-1}(s) \right)$.

Step $n+1-i$, where $i \in \{2,...,n-1\}$ : Find $k_i \in [x_i^*, a_i)$ such that $\min_{k_i \leq s \leq a_i} \left( \min(a_{i+1} - k_{i+1}, f_i(s)) - \min(a_i - s, F_{i-1}) \right) > 0$ and $\max_{0 \leq s \leq k_i} \left( s - \min(a_{i+1} - k_{i+1}, f_i(s)) + \min(a_i - s, F_{i-1}) \right) \leq k_i$, where $F_{i-1} := \max_{s \in [0, a_{i-1}]} \left( f_{i-1}(s) \right)$.

Step $n$: Find $c_1 \in [x_1^*, a_1)$ such that $\min_{c_1 \leq s \leq a_1} \left( \min(a_2 - k_2, f_1(s)) - \min(a_1 - s, v) \right) > 0$ and $\max_{0 \leq s \leq c_1} \left( s - \min(a_2 - k_2, f_1(s)) + \min(a_1 - s, v) \right) \leq c_1$.

Step $n+i-1$, where $i \in \{2,...,n-1\}$ : Find $c_i \in [x_i^*, k_i]$ such that $\min_{c_i \leq s \leq k_i} \left( \min(a_{i+1} - k_{i+1}, f_i(s)) - \min(a_i - s, F_{i-1}) \right) > 0$ and $\max_{0 \leq s \leq c_i} \left( s - \min(a_{i+1} - k_{i+1}, f_i(s)) + \min(a_i - s, F_{i-1}) \right) \leq c_i$, where $F_{i-1} := \max_{s \in [0, c_{i-1}]} \left( f_{i-1}(s) \right)$.

Step $2n-1$: Find $c_n \in [x_n^*, k_n]$ such that $\min_{c_n \leq s \leq k_n} \left( f_n(s) - \min(a_n - s, F_{n-1}) \right) > 0$ and $\max_{0 \leq s \leq c_n} \left( s - f_n(s) + \min(a_n - s, F_{n-1}) \right) \leq c_n$, where $F_{n-1} := \max_{s \in [0, c_{n-1}]} \left( f_{n-1}(s) \right)$.

*Furthermore, assume that there exist constants* $\lambda_i \in [0,1)$ *(* $i = 1,...,n$ *) such that*

$$\left| s - x_i^* - \min(a_{i+1} - c_{i+1}, f_i(s)) + \min(a_i - s, v) \right| \leq \lambda_i \left| s - x_i^* \right|, \text{ for } s \in [0, c_i],\ i = 1,...,n-1 \tag{4.26}$$

$$\left| s - x_i^* - \min(a_{i+1}, f_i(s)) + \min(a_i - s, v) \right| \leq \lambda_i \left| s - x_i^* \right|, \text{ for } s \in [0, c_i],\ i = 1,...,n-1 \tag{4.27}$$

$$\left| s - x_n^* - f_n(s) + \min(a_n - s, v) \right| \leq \lambda_n \left| s - x_n^* \right|, \text{ for all } s \in [0, c_n] \tag{4.28}$$

*Moreover, assume that there exists a constant* $L \geq 0$ *such that inequalities (4.7) hold. Then the equilibrium point* $x^* = (x_1^*,...,x_n^*) \in [0, a_1] \times \ldots [0, a_n]$ *is GES for the disturbance-free system (4.4), (4.5), (4.6).*

All steps of the algorithm of Corollary 4.4 can be performed due to the fact that the functions $f_i : [0, a_i] \to \Re_+$ are continuous functions with $0 < f_i(s) < a_{i+1}$ for all $s \in (0, a_i]$ ($i = 1,...,n-1$). Moreover, all steps of the algorithm of Corollary 4.4 can be easily performed numerically. Once we have found all constants $c_i$ ($i = 1,...,n$), it is straightforward to determine numerically the values of the constants $\lambda_i$ ($i = 1,...,n$) for which inequalities (4.26), (4.27), (4.28) hold: if all constants $\lambda_i$ ($i = 1,...,n$) are less than 1, then we can guarantee GES for the equilibrium point $x^* = (x_1^*,...,x_n^*) \in [0, a_1] \times \ldots [0, a_n]$ of the disturbance-free system (4.4), (4.5), (4.6).

**Example 4.5:** We consider the network (4.4), (4.5), (4.6) with

$$n = 5,\ a_i = 10\ (i = 1,...,5),\ f_i(s) = f(s) := \begin{cases} 0.5s & \text{for } s \in [0,5] \\ 3 - 0.1s & \text{for } s \in (5,10] \end{cases}\ (i = 1,...,4),$$

$$f_5(s) := \begin{cases} 0.4s & \text{for } s \in [0,5] \\ 2 - p(s-5) & \text{for } s \in (5,10] \end{cases},\ v = 1 \tag{4.1}$$

where $p \in [0, 0.4)$. For this network we have

$$x_i^* = 2\ (i = 1,...,4),\ x_5^* = 2.5 \tag{4.2}$$

and assumption (H) holds for the equilibrium point $x^* \in \Re^5$. We are using Corollary 4.4 in order to answer the following question:

"For what values of $p \in [0,0.4)$ the equilibrium point $x^* \in \Re^5$ is GES?"

The algorithm of Corollary 4.4 was performed numerically for various values of $p \in [0,0.4)$ in the following way: for a given integer $N > 0$ a grid of points $s_i = i\frac{a}{N}$ ($i = 0,1,...,N$) was generated. Then $c_i$ ($i = 1,...,5$) and $k_i$ ($i = 2,...,5$) were chosen to be equal to the smallest grid point $s_j = j\frac{a}{N}$ so that $\min_{j \le l \le N}\left( q\left( l\frac{a}{N} \right) \right) > 0$ and $\max_{0 \le l \le j}\left( l\frac{a}{N} - q\left( l\frac{a}{N} \right) \right) < j\frac{a}{N}$, where:

- $q(s) := f_5(s) - \min(10 - s, 2.5)$ and $F := 2.5$ for the determination of $k_5$,
- $q(s) := \min(10 - k_{i+1}, f(s)) - \min(10 - s, 2.5)$, for the determination of $k_i$ ($i = 4,3,2$),
- $q(s) := \min(10 - k_2, f(s)) - \min(10 - s, 1)$ for the determination of $c_1$,
- $q(s) := \min(10 - k_{i+1}, f(s)) - \min(10 - s, F_{i-1})$, $F_{i-1} := \max_{s \in [0, c_{i-1}]} (f(s))$ for the determination of $c_i$ ($i = 2,3,4$),
- $q(s) := f_5(s) - \min(10 - s, F_4)$ and $F_4 := \max_{s \in [0, c_4]} (f(s))$, for the determination of $c_5$.

For various values of $N > 0$, it was found that there exists $p_N > 0$ such that the assumptions of Corollary 4.4 hold with $\lambda_i = 0.5$ ($i = 1,...,4$) and $\lambda_5 = 0.6$ for all $p \in [0, p_N]$. More specifically, we obtained:

- $p_{100} = 0.18918$,
- $p_{1000} = 0.244332$,
- $p_{2000} = 0.247176$,

indicating a sequence that tends to 0.25 as $N \to +\infty$. Simulations show that the equilibrium point $x^* \in \Re^5$ is GES for all $p \in [0,0.25)$. Indeed, for $p = 0.25$ there exist additional (a continuum of) equilibrium points to which many solutions are attracted and therefore, the equilibrium point $x^* \in \Re^5$ cannot be GES for $p = 0.25$.

The results of this example show that the sufficient conditions of Corollary 4.4 are virtually exact in this case. However, in general, it is expected to have a trade-off between the goal of obtaining easily checkable sufficient conditions for GES and the reduction of conservatism. ◁

## 5. Concluding Remarks

Sufficient conditions for GAS and GES have been given, by means of vector Lyapunov functions. The conditions can be applied to nonlinear, large-scale, uncertain discrete-time systems. The obtained results were applied to traffic networks for the derivation of sufficient conditions of GAS of the uncongested equilibrium point of the network. Specific results and algorithms were provided for freeway models. Various examples illustrated the applicability of the obtained results.

The results of the present paper can be used for different purposes for future research:

- for the derivation of feedback laws which stabilize the uncongested equilibrium point,
- for the study of the dynamic behavior of traffic networks under the effect of external disturbances (varying inflows),
- for the study of more complicated freeway models divided in $n \ge 3$ cells, each with one on-ramp and one off-ramp.

All the above topics can give important results for transportation and traffic theory.

**Acknowledgments**

The research leading to these results has received funding from the European Research Council under the European Union's Seventh Framework Programme (FP/2007-2013) / ERC Grant Agreement n. [321132], project TRAMAN21.